\documentclass[12pt,oneside,reqno]{amsart}
\usepackage{txfonts}
\usepackage{bbm}
\usepackage{amsmath}
\usepackage{graphicx}
\usepackage{mathrsfs}
\usepackage{stmaryrd}
\usepackage{amsfonts}
\usepackage{enumerate,amsmath,amssymb,amsthm}

\numberwithin{equation}{section}

\newcommand{\be}{\begin{eqnarray}}
\newcommand{\ee}{\end{eqnarray}}
\newcommand{\ce}{\begin{eqnarray*}}
\newcommand{\de}{\end{eqnarray*}}
\newtheorem{theorem}{Theorem}[section]
\newtheorem{lemma}[theorem]{Lemma}
\newtheorem{remark}[theorem]{Remark}
\newtheorem{definition}[theorem]{Definition}
\newtheorem{proposition}[theorem]{Proposition}
\newtheorem{Examples}[theorem]{Example}
\newtheorem{corollary}[theorem]{Corollary}

\def\u{{\mathbf{u}}}
\def\eps{\varepsilon}

\def\a{\alpha}

\def\p{\partial}

\def\[{{\Big[}}
\def\]{{\Big]}}
\def\<{{\langle}}
\def\>{{\rangle}}
\def\({{\Big(}}
\def\){{\Big)}}

\def\bx{{\mathbf{x}}}

\def\dif{{\mathord{{\rm d}}}}

\def\no{\nonumber}
\def\={&\!\!=\!\!&}
\def\bt{\begin{theorem}}
\def\et{\end{theorem}}
\def\bl{\begin{lemma}}
\def\el{\end{lemma}}
\def\br{\begin{remark}}
\def\er{\end{remark}}

\def\bd{\begin{definition}}
\def\ed{\end{definition}}
\def\bp{\begin{proposition}}
\def\ep{\end{proposition}}
\def\bc{\begin{corollary}}
\def\ec{\end{corollary}}
\def\bx{\begin{Examples}}
\def\ex{\end{Examples}}

\def\cH{{\mathcal H}}

\def\cJ{{\mathcal J}}

\def\cM{{\mathcal M}}

\def\cR{{\mathcal R}}

\def\mE{{\mathbb E}}

\def\mH{{\mathbb H}}
\def\mI{{\mathbb I}}

\def\mL{{\mathbb L}}

\def\mN{{\mathbb N}}

\def\mR{{\mathbb R}}

\def\sF{{\mathscr F}}

\def\sS{{\mathscr S}}

\def\geq{\geqslant}
\def\leq{\leqslant}

\allowdisplaybreaks

\begin{document}

\title{Stochastic Homeomorphism Flows of SDEs with Singular Drifts and Sobolev Diffusion Coefficients}

\date{}
\author{Xicheng Zhang }

\dedicatory{
School of Mathematics and Statistics,
Wuhan University, Wuhan, Hubei 430072, P.R.China\\
Email: XichengZhang@gmail.com
 }

\begin{abstract}
In this paper we prove the stochastic homeomorphism flow property and the
strong Feller property for stochastic differential equations with sigular time dependent drifts and
Sobolev diffusion coefficients. Moreover, the local well posedness under local assumptions
are also obtained. In particular, we extend Krylov and R\"ockner's results in \cite{Kr-Ro}
to the case of non-constant diffusion coefficients.
\end{abstract}

\thanks{{\it AMS Classification(2000)}: 60H15}

\maketitle \rm

\section{Introduction and Main Result}

Consider the following stochastic differential equation (SDE) in $\mR^d$:
\begin{align}
\dif X_t=b_t(X_t)\dif t+\sigma_t(X_t)\dif W_t,\label{SDE}
\end{align}
where $b:\mR_+\times\mR^d\to\mR^d$ and $\sigma:\mR_+\times\mR^d\to\mR^d\times\mR^d$ are two Borel measurable functions, and
$\{W_t\}_{t\geq 0}$
is a $d$-dimensional standard Brownian motion defined on some complete filtered probability space $(\Omega,\sF,P;(\sF_t)_{t \geq 0})$.
When $\sigma$ is Lipschitz continuous in $x$ uniformly with respect to $t$ and $b$ is bounded measurable, Veretennikov \cite{Ve} first
proved the existence of a unique strong solution for SDE (\ref{SDE}).
Recently, Krylov and R\"ockner \cite{Kr-Ro} proved the existence and uniqueness of strong solutions for
SDE (\ref{SDE}) with  $\sigma\equiv\mI_{d\times d}$ and
\begin{align}
\int^T_0\left(\int_{\mR^d}|b_t(x)|^p\dif x\right)^{\frac{q}{p}}\dif t<+\infty, \ \ \forall T>0,\label{BP3}
\end{align}
provided that
\begin{align}
\frac{d}{p}+\frac{2}{q}<1.\label{BP4}
\end{align}
More recently, following \cite{Kr-Ro}, Fedrizzi and Flandoli \cite{Fe-Fl} proved the $\alpha$-H\"older continuity of $x\mapsto X_t(x)$ for any $\alpha\in(0,1)$
basing on Girsanov's theorem and Khasminskii's estimate.
In the case of non-constant and non-degenerate diffusion coefficient,
the present author \cite{Zh0} proved the pathwise uniqueness for SDE (\ref{SDE})
under stronger integrability assumptions on $b$ and $\sigma$ (see also \cite{Gy-Ma} for Lipschitz $\sigma$ and unbounded $b$).
Moreover, there are many works recently devoted to the study of stochastic homeomorphism (or diffeomorphism)
flow property of SDE (\ref{SDE}) under various non-Lipschitz assumptions on coefficients (see \cite{Fa-Im-Zh,Zh2,Fl-Gu-Pr} and references therein).

We first introduce the class of local strong solutions for SDE (\ref{SDE}). Let $\tau$ be any ($\sF_t$)-stopping time
and $\xi$ any $\sF_0$-measurable $\mR^d$-valued random variable. Let $\sS^\tau_{b,\sigma}(\xi)$
be the class of all $\mR^d$-valued ($\sF_t$)-adapted continuous stochastic process $X_t$ on $[0,\tau)$ satisfying
$$
P\left\{\omega: \int^T_0|b_s(X_s(\omega))|\dif s+\int^T_0|\sigma_s(X_s(\omega))|^2\dif s<+\infty,
\forall T\in[0,\tau(\omega))\right\}=1,
$$
and such that
$$
X_t=\xi+\int^t_0b_s(X_s)\dif s+\int^t_0\sigma_s(X_s)\dif W_s,\ \ \forall t\in[0,\tau),\ \ a.s.
$$
We now state our main result as follows:
\bt\label{Main}
In addition to (\ref{BP3}) with $p,q\in(1,\infty)$ satisfying (\ref{BP4}), we also assume that
\begin{enumerate}[{\bf (H$^\sigma_1$)}]
\item $\sigma_t(x)$ is uniformly continuous in $x\in\mR^d$ locally uniformly with respect to $t\in\mR_+$,
and there exist positive constants $K$ and $\delta$ such that for all $(t,x)\in\mR_+\times\mR^d$,
$$
\delta|\lambda|^2\leq\sum_{ik}|\sigma^{ik}_t(x)\lambda^i|^2\leq K|\lambda|^2,\ \ \forall\lambda\in\mR^d;
$$
\item $|\nabla\sigma_t|\in L^{q}_{loc}(\mR_+;L^{p}(\mR^d))$ with the same $p,q$ as required on $b$, where $\nabla$ denotes the
generalized gradient with respect to $x$.
\end{enumerate}
Then for any ($\sF_t$)-stopping time $\tau$ (possibly being infinity)
and $x\in\mR^d$, there exists a unique strong solution $X_t(x)\in\sS^\tau_{b,\sigma}(x)$
to SDE (\ref{SDE}), which means that for any $X_t(x),Y_t(x)\in\sS^\tau_{b,\sigma}(x)$,
$$
P\{\omega: X_t(\omega,x)=Y_t(\omega,x), \forall t\in[0,\tau(\omega))\}=1.
$$
Moreover, for almost all $\omega$ and all $t\geq 0$,
$$
x\mapsto X_t(\omega,x)\mbox{ is a homeomorphism on $\mR^d$},
$$
and for any $t>0$ and bounded measurable function $\phi$,  $x,y\in\mR^d$,
$$
|\mE\phi(X_t(x))-\mE\phi(X_t(y))|\leq C_t\|\phi\|_\infty|x-y|,
$$
where $C_t>0$ satisfies $\lim_{t\to 0}C_t=+\infty$.
\et
\br
The uniqueness proven in this theorem means local uniqueness. We want to emphasize that global uniqueness
can not imply local uniqueness since local solution can not in general be extended to a global solution.
\er
By localization technique (cf. \cite{Zh0}), as a corollary of Theorem \ref{Main},
we have the following existence and uniqueness of local strong solutions.
\bt
Assume that for any $n\in\mN$ and some $p_n,q_n\in(1,\infty)$ satisfying (\ref{BP4}),
\begin{enumerate}[(i)]
\item $|b_t|, |\nabla\sigma_t|\in L^{q_n}_{loc}(\mR_+;L^{p_n}(B_n))$, where $B_n:=\{x\in\mR^d: |x|\leq n\}$;
\item $\sigma^{ik}_t(x)$ is uniformly continuous in $x\in B_n$ uniformly with respect to $t\in[0,n]$,
and there exist positive constants $\delta_n$ such that for all $(t,x)\in[0,n]\times B_n$,
$$
\sum_{ik}|\sigma^{ik}_t(x)\lambda^i|^2\geq\delta_n|\lambda|^2,\ \ \forall\lambda\in\mR^d.
$$
\end{enumerate}
Then for any $x\in\mR^d$, there exist an ($\sF_t$)-stopping time $\zeta(x)$ (called explosion time) and a unique strong solution
$X_t(x)\in\sS^{\zeta(x)}_{b,\sigma}(x)$ to SDE (\ref{SDE}) such that on $\{\omega: \zeta(\omega,x)<+\infty\}$,
\begin{align}
\lim_{t\uparrow\zeta(x)}X_t(x)=+\infty,\ \ a.s.\label{HH1}
\end{align}
\et
\begin{proof}
For each $n\in\mN$, let $\chi_n(t,x)\in[0,1]$ be a nonnegative smooth function in $\mR_+\times\mR^d$ 
with $\chi_n(t,x)=1$ for all $(t,x)\in [0,n]\times B_n$
and $\chi_n(t,x)=0$ for all $(t,x)\notin [0,n+1]\times B_{n+1}$. Let
$$
b^n_t(x):=\chi_n(t,x)b_t(x)
$$
and
$$
\sigma^n_t(x):=\chi_{n+1}(t,x)\sigma_t(x)
+(1-\chi_n(t,x))\left(1+\sup_{(t,x)\in [0,n+2]\times B_{n+2}}|\sigma_t(x)|\right)\mI_{d\times d}.
$$
By Theorem \ref{Main}, for each $x\in\mR^d$, there exists a unique strong solution $X^n_t(x)\in\sS^\infty_{b^n,\sigma^n}(x)$
to SDE (\ref{SDE}) with coefficients $b^n$ and $\sigma^n$. For $n\geq k$, define
$$
\tau_{n,k}(x,\omega):=\inf\{t\geq 0: |X^n_t(\omega,x)|\geq k\}\wedge n.
$$
It is easy to see that
$$
X^n_t(x), X^k_t(x)\in\sS^{\tau_{n,k}(x)}_{b^k,\sigma^k}(x).
$$
By the local uniqueness proven in Theorem \ref{Main}, we have
$$
P\{\omega:X^n_t(\omega,x)= X^k_t(\omega,x), \forall t\in[0,\tau_{n,k}(x,\omega))\}=1,
$$
which implies that for $n\geq k$,
$$
\tau_{k,k}(x)\leq\tau_{n,k}(x)\leq\tau_{n,n}(x), \ \ a.s.
$$
Hence, if we let $\zeta_k(x):=\tau_{k,k}(x)$, then $\zeta_k(x)$ is an increasing sequence of ($\sF_t$)-stopping times and
for $n\geq k$,
$$
P\{\omega: X^n_t(x,\omega)=X^k_t(x,\omega),\ \ \forall t\in[0,\zeta_k(x,\omega))\}=1.
$$
Now, for each $k\in\mN$, we can define $X_t(x,\omega)=X^k_t(x,\omega)$ for $t<\zeta_k(x,\omega)$ and $\zeta(x)=\lim_{k\to\infty}\zeta_k(x)$.
It is clear that $X_t(x)\in\sS^{\zeta(x)}_{b,\sigma}(x)$ and (\ref{HH1}) holds.
\end{proof}

The aim of this paper is now to prove Theorem \ref{Main}. We organize it as follows: In Section 2,
we prove two new estimates of Krylov's type, which is the key point for our proof and has some independent interest.
In Section 3, we prove Theorem \ref{Main} in the case of $b=0$. For the stochastic homeomorphism flow, we adopt
Kunita's simple argument (cf. \cite{Ku}).
For the strong Feller property, we use Bismut-Elworthy-Li's formula (cf. \cite{El-Li}).
In Section 4, we use Zvonkin's transformation
to fully prove Theorem \ref{Main}. In Appendix, we recall some well known facts used in the present paper.

\section{Two estimates of Krylov's type }

We first introduce some spaces and notations. For $p,q\in[1,\infty)$ and $0\leq S<T<\infty$, we denote by $\mL^q_p(S,T)$ the space of
all real Borel measurable functions on $[S,T]\times\mR^d$ with the norm
$$
\|f\|_{\mL^q_p(S,T)}:=\left(\int^T_S\left(\int_{\mR^d}f(t,x)^p\dif x\right)^{\frac{q}{p}}\right)^{\frac{1}{q}}<+\infty.
$$
For $m\in\mN$ and $p\geq 1$, let $H^m_p$ be the usual Sobolev space over $\mR^d$ with the norm
$$
\|f\|_{H^m_p}:=\sum_{k=0}^m\|\nabla^k f\|_{L^p}<+\infty,
$$
where $\nabla$ denotes the gradient operator, and $\|\cdot\|_{L^p}$ is the usual $L^p$-norm.
We also introduce for $0\leq S<T<\infty$,
$$
\mH^{2,q}_p(S,T)=L^q(S,T;H^2_p),
$$
and the space $\cH^{2,q}_p(S,T)$ consisting of function $u=u(t)$ defined on $[S,T]$ with values in the space
of distributions on $\mR^d$ such that $u\in\mH^{2,q}_p(S,T)$ and $\p_t u\in\mL^q_p(S,T)$.
For simplicity, we write
$$
\mL^q_p(T)=\mL^q_p(0,T),\ \ \ \mH^{2,q}_p(T)=\mH^{2,q}_p(0,T),\ \ \ \cH^{2,q}_p(T)=\cH^{2,q}_p(0,T)
$$
and
\begin{align}
L_tu(x):=\tfrac{1}{2}\sigma^{ik}_t(x)\sigma^{jk}_t(x)\p_i\p_ju(x)+b^i_t(x)\p_iu(x).\label{LL2}
\end{align}
Here and below, we use the convention that the repeated
indices in a product will be summed automatically. Moreover, the letter $C$ will denote an unimportant constant, whose dependence on the
functions or parameters can be traced from the context.

We first prove the following estimate of Krylov's type (cf. \cite[p.54, Theorem 4]{Kr}).
\bt\label{Krylov}
Suppose that $\sigma$ satisfies {\bf (H$^\sigma_1$)} and $b$ is bounded measurable. Fix an ($\sF_t$)-stopping time $\tau$ and an
$\sF_0$-measurable $\mR^d$-valued random variable $\xi$ and let $X_t\in\sS^\tau_{b,\sigma}(\xi)$.
Given $T_0>0$ and $p,q\in(1,\infty)$ with
\begin{align}
\frac{d}{p}+\frac{2}{q}<2,\label{BB6}
\end{align}
there exists a positive constant $C=C(K,\delta,d,p,q, T_0, \|b\|_\infty)$
such that for all $f\in\mL^q_p(T_0)$ and $0\leq S<T\leq T_0$,
\begin{align}
\mE\left(\int^{T\wedge\tau}_{S\wedge\tau} f(s,X_s)\dif s\Bigg|_{\sF_S}\right)\leq C\|f\|_{\mL^q_p(S,T)}.\label{Es1}
\end{align}
\et
\begin{proof}
Let $r=d+1$. Since $\mL^r_r(T_0)\cap\mL^q_p(T_0)$ is dense in $\mL^q_p(T_0)$, it suffices to prove (\ref{Es1}) for
$$
f\in\mL^r_r(T_0)\cap\mL^q_p(T_0).
$$
Fix $T\in[0,T_0]$. By Theorem \ref{Th1} in appendix, there exists a unique solution $u\in\cH^{2,r}_r(T)\cap\cH^{2,q}_p(T)$
for the following backward PDE on $[0,T]$:
$$
\p_t u(t,x)+L_tu(t,x)=f(t,x),\ \ \ u(T,x)=0.
$$
Moreover, for some constant $C=C(K,\delta,d,p,q, T_0, \|b\|_\infty)$,
\begin{align}
\|\p_tu\|_{\mL^r_r(S,T)}+\|u\|_{\mH^{2,r}_r(S,T)}\leq C\|f\|_{\mL^r_r(S,T)}, \ \ \forall S\in[0,T]\label{BB7}
\end{align}
and
$$
\|\p_tu\|_{\mL^q_p(S,T)}+\|u\|_{\mH^{2,q}_p(S,T)}\leq C\|f\|_{\mL^q_p(S,T)}, \ \ \forall S\in[0,T].
$$
In particular, by (\ref{BB6}) and \cite[Lemma 10.2]{Kr-Ro},
\begin{align}
\sup_{(t,x)\in[S,T]\times\mR^d}|u(t,x)|\leq C\|f\|_{\mL^q_p(S,T)}.\label{BP7}
\end{align}

Let $\rho$ be a nonnegative smooth function in $\mR^{d+1}$ with support in $\{x\in\mR^{d+1}: |x|\leq 1\}$ and
$\int_{\mR^{d+1}}\rho(t,x)\dif t\dif x=1$.
Set $\rho_n(t,x):=n^{d+1}\rho(nt,nx)$ and extend $u(s)$ to $\mR$ by setting
$u(s,\cdot)=0$ for $s\geq T$ and $u(s,\cdot)=u(0,\cdot)$ for $s\leq 0$. Define
\begin{align}
u_n(t,x):=\int_{\mR^{d+1}}u(s,y)\rho_n(t-s,x-y)\dif s\dif y\label{BP6}
\end{align}
and
$$
f_n(t,x):=\p_t u_n(t,x)+L_tu_n(t,x).
$$
Then by (\ref{BB7}) and the property of convolutions, we have
\begin{align*}
\|f_n-f\|_{\mL^r_r(T)}&\leq\|\p_t(u_n-u)\|_{\mL^r_r(T)}+\|b^i\p_i(u_n-u)\|_{\mL^r_r(T)}
+K\|\p_i\p_j(u_n-u)\|_{\mL^r_r(T)}\no\\
&\leq\|\p_t(u_n-u)\|_{\mL^r_r(T)}+\|b\|_\infty\|\nabla(u_n-u)\|_{\mL^r_r(T)}
+K\|u_n-u\|_{\mH^{2,r}_r(T)}\no\\
&\leq\|\p_t(u_n-u)\|_{\mL^r_r(T)}+C\|u_n-u\|_{\mH^{2,r}_r(T)}\rightarrow 0\mbox{ as $n\to\infty$}.
\end{align*}
So, by the classical Krylov's estimate (cf. \cite[Lemma 5.1]{Kr4} or \cite[Lemma 3.1]{Gy-Ma}), we have
\begin{align}
\lim_{n\to\infty}\mE\left(\int^{T\wedge\tau}_0 |f_n(s,X_s)-f(s,X_s)|\dif s\right)
\leq\lim_{n\to\infty}\|f_n-f\|_{\mL^r_r(T)}=0.\label{BP5}
\end{align}

Now using It\^o's formula for $u_n(t,x)$, we have
$$
u_n(t,X_t)=u_n(0,X_0)+\int^t_0 f_n(s,X_s)\dif s+\int^t_0\p_i u_n(s,X_s)\sigma^{ik}_s(X_s)\dif W^k_s,\ \ \forall t<\tau.
$$
In view of
$$
\sup_{s,x}|\p_i u_n(s,x)|\leq C_n,
$$
by Doob's optional theorem, we have
$$
\mE\left[\int^{T\wedge\tau}_{S\wedge\tau}\p_i u_n(s,X_s)\sigma^{ik}_s(X_s)\dif W^k_s\Bigg|_{\sF_S}\right]=0.
$$
Hence,
\begin{align}
\mE\left(\int^{T\wedge\tau}_{S\wedge\tau} f_n(s,X_s)\dif s\Bigg|_{\sF_S}\right)
&=\mE\Bigg[(u_n(T\wedge\tau,X_{T\wedge\tau})-u_n(S\wedge\tau,X_{S\wedge\tau}))\Big|_{\sF_S}\Bigg]\label{LL1}\\
&\leq 2\sup_{(t,x)\in[S,T]\times\mR^d}|u_n(t,x)|\leq 2\sup_{(t,x)\in[S,T]\times\mR^d}|u(t,x)|
\stackrel{(\ref{BP7})}{\leq} C\|f\|_{\mL^q_p(S,T)}\no.
\end{align}
The proof is thus completed by (\ref{BP5}) and letting $n\to\infty$.
\end{proof}

Next, we want to relax the boundedness assumption on $b$. The price to pay is that a stronger integrability assumption is required.
\bt\label{Krylov1}
Suppose that $\sigma$ satisfies {\bf (H$^\sigma_1$)} and $b\in L^q(\mR_+,L^p(\mR^d))$ provided with
\begin{align}
\frac{d}{p}+\frac{2}{q}<1.\label{BB8}
\end{align}
 Fix an ($\sF_t$)-stopping time $\tau$ and an
$\sF_0$-measurable $\mR^d$-valued random variable $\xi$ and let $X_t\in\sS^\tau_{b,\sigma}(\xi)$.
Given $T_0>0$, there exists a positive constant $C=C(K,\delta,d,p,q, T_0, \|b\|_{\mL^q_p(T_0)})$
such that for all $f\in\mL^q_p(T_0)$ and $0\leq S<T\leq T_0$,
\begin{align}
\mE\left(\int^{T\wedge\tau}_{S\wedge\tau} f(s,X_s)\dif s\Bigg|_{\sF_S}\right)\leq C\|f\|_{\mL^q_p(S,T)}. \label{BB9}
\end{align}
\et
\begin{proof}
Following the proof of Theorem \ref{Krylov}, we let $r=d+1$ and
assume that
$$
f\in\mL^r_r(T_0)\cap\mL^q_p(T_0).
$$
Below, for $N>0$, we write
$$
L^N_tu(x):=\tfrac{1}{2}\sigma^{ik}_t(x)\sigma^{jk}_t(x)\p_i\p_ju(x)+1_{\{|b_t(x)|\leq N\}}b^i_t(x)\p_iu(x).
$$
Fix $T\in[0,T_0]$. By Theorem \ref{Th1}, there exists a unique solution $u\in\cH^{2,r}_r(T)\cap\cH^{2,q}_p(T)$
for the following backward PDE on $[0,T]$:
$$
\p_t u(t,x)+L^N_tu(t,x)=f(t,x),\ \ \ u(T,x)=0.
$$
Moreover, for some constant $C_1=C_1(K,\delta,d,p,q, T_0, N)$,
\begin{align}
\|\p_tu\|_{\mL^r_r(S,T)}+\|u\|_{\mH^{2,r}_r(S,T)}\leq C_1\|f\|_{\mL^r_r(S,T)}, \ \ \forall S\in[0,T],\label{BB5}
\end{align}
and for some constant $C_2=C_2(K,\delta,d,p,q, T_0, \|b\|_{\mL^q_p(T)})$,
$$
\|\p_tu\|_{\mL^q_p(S,T)}+\|u\|_{\mH^{2,q}_p(S,T)}\leq C_2\|f\|_{\mL^q_p(S,T)}, \ \ \forall S\in[0,T].
$$
In particular, by (\ref{BB8}) and \cite[Lemma 10.2]{Kr-Ro},
\begin{align}
\sup_{(t,x)\in[S,T]\times\mR^d}|u(t,x)|
+\sup_{(t,x)\in[S,T]\times\mR^d}|\nabla u(t,x)|\leq C_2\|f\|_{\mL^q_p(S,T)}.\label{BP77}
\end{align}

For $R>0$, define
$$
\tau_R:=\inf\left\{t\in[0,\tau): \int^t_0|b_s(X_s)|\dif s\geq R\right\}.
$$
Let $u_n$ be defined by (\ref{BP6}).
As in the proof of Theorem \ref{Krylov} (see (\ref{LL1})), by (\ref{BP77}), we have
\begin{align}
\mE\left(\int^{T\wedge\tau_R}_{S\wedge\tau_R} (\p_su_n+L_su_n)(s,X_s)\dif s\Bigg|_{\sF_S}\right)
\leq C_2\|f\|_{\mL^q_p(S,T)}.\label{BB4}
\end{align}
Now if we set
$$
f^N_n(t,x):=\p_t u_n(t,x)+L^N_tu_n(t,x),
$$
then
\begin{align*}
\mE\left(\int^{T\wedge\tau_R}_{S\wedge\tau_R} f^N_n(s,X_s)\dif s\Bigg|_{\sF_S}\right)
&=\mE\left(\int^{T\wedge\tau_R}_{S\wedge\tau_R} (\p_su_n+L_su_n)(s,X_s)\dif s\Bigg|_{\sF_S}\right)\\
&\quad-\mE\left(\int^{T\wedge\tau_R}_{S\wedge\tau_R} 1_{\{|b_s(X_s)|>N\}}b^i_s(X_s)\p_iu_n(s,X_s)\dif s\Bigg|_{\sF_S}\right).
\end{align*}
Hence, by (\ref{BP77}) and (\ref{BB4}),
\begin{align}
\mE\left(\int^{T\wedge\tau_R}_{S\wedge\tau_R} f^N_n(s,X_s)\dif s\Bigg|_{\sF_S}\right)\leq
C\|f\|_{\mL^q_p(S,T)}+C\mE\left(\int^{T\wedge\tau_R}_{S\wedge\tau_R}
1_{\{|b_s(X_s)|>N\}}|b_s(X_s)|\dif s\Bigg|_{\sF_S}\right),\label{BB2}
\end{align}
where $C=C(K,\delta,d,p,q, T_0, \|b\|_{\mL^q_p(T_0)})$ is independent of $n$ and $R,N$.
Observe that for fixed $N>0$, by (\ref{BB5}),
$$
\lim_{n\to\infty}\|f^N_n-f\|_{\mL^r_r(T)}=0,
$$
and for fixed $R>0$, by the dominated convergence theorem,
$$
\lim_{N\to\infty}\mE\left(\int^{T\wedge\tau_R}_{S\wedge\tau_R}
1_{\{|b_s(X_s)|>N\}}|b_s(X_s)|\dif s\right)=0.
$$
Taking limits  for both sides of (\ref{BB2}) in order: $n\to\infty$, $N\to\infty$ and $R\to\infty$, we obtain (\ref{BB9}).
\end{proof}

\section{SDE with Sobolev diffusion coefficient and zero drift}

In this section we consider the following SDE without drift:
\begin{align}
X_t(x)=x+\int^t_0\sigma_s(X_s(x))\dif W_s.\label{DSDE}
\end{align}
We first prove that:
\bt\label{Path}
Under {\bf (H$^\sigma_1$)} and {\bf (H$^\sigma_2$)}, the local pathwise uniqueness holds for SDE (\ref{DSDE}). More precisely,
for any ($\sF_t$)-stopping time $\tau$ (possibly being infinity) and $x\in\mR^d$, let $X_t, Y_t\in\sS^\tau_{0,\sigma}(x)$, then
$$
P\{\omega: X_t(\omega)=Y_t(\omega),\forall t\in[0,\tau(\omega))\}=1.
$$
In particular, there exists a unique strong solution for SDE (\ref{DSDE}).
\et
\begin{proof}
Set $Z_t:=X_t-Y_t$. By It\^o's formula, we have
\begin{align*}
|Z_{t\wedge\tau}|^2=2\int^{t\wedge\tau}_0\<Z_s,[\sigma_s(X_s)-\sigma_s(Y_s)]\dif W_s\>
+\int^{t\wedge\tau}_0\|\sigma_s(X_s)-\sigma_s(Y_s)\|^2\dif s.
\end{align*}
If we set
$$
M_t:=2\int^t_0\frac{\<Z_s,[\sigma_s(X_s)-\sigma_s(Y_s)]\dif W_s\>}{|Z_s|^2}
$$
and
$$
A_t:=\int^t_0\frac{\|\sigma_s(X_s)-\sigma_s(Y_s)\|^2}{|Z_s|^2}\dif s,
$$
then
$$
|Z_{t\wedge\tau}|^2=\int^{t\wedge\tau}_0|Z_s|^2\dif (M_s+ A_s).
$$
Here and below, we use the convention that $\frac{0}{0}\equiv0$.
Thus, if we can show that $t\mapsto M_{t\wedge\tau}+A_{t\wedge\tau}$ is a continuous semimartingale,
then the uniqueness follows. For this, it suffices to prove that for any $t\geq 0$,
$$
\mE|M_{t\wedge\tau}|^2<+\infty, \ \ \mE A_{t\wedge\tau}<+\infty.
$$
Set
$$
\sigma^n_s(x):=\sigma_s*\rho_n(x),
$$
where $\rho_n$ is a mollifier in $\mR^d$ as used in Theorem \ref{Krylov}.
By Fatou's lemma, we have
\begin{align*}
\mE A_{t\wedge\tau}&\leq\varliminf_{\eps\downarrow 0}\mE\int^{t\wedge\tau}_0\frac{\|\sigma_s(X_s)-\sigma_s(Y_s)\|^2}{|Z_s|^2}\cdot1_{|Z_s|>\eps}\dif s\\
&\leq 3\Bigg(\varliminf_{\eps\downarrow 0}\sup_{n\in\mN}\mE\int^{t\wedge\tau}_0
\frac{\|\sigma^n_s(X_s)-\sigma^n_s(Y_s)\|^2}{|Z_s|^2}\cdot1_{|Z_s|>\eps}\dif s\\
&\quad+\varliminf_{\eps\downarrow 0}\lim_{n\to\infty}\mE\int^{t\wedge\tau}_0
\frac{\|\sigma^n_s(X_s)-\sigma_s(X_s)\|^2}{|Z_s|^2}\cdot1_{|Z_s|>\eps}\dif s\\
&\quad+\varliminf_{\eps\downarrow 0}\lim_{n\to\infty}\mE\int^{t\wedge\tau}_0
\frac{\|\sigma^n_s(Y_s)-\sigma_s(Y_s)\|^2}{|Z_s|^2}\cdot1_{|Z_s|>\eps}\dif s\Bigg)\\
&=:3(I_1(t)+I_2(t)+I_3(t)).
\end{align*}
By estimate (\ref{Es1}), we have
\begin{align*}
I_2(t)&\leq\varliminf_{\eps\downarrow 0}\frac{1}{\eps^2}\lim_{n\to\infty}\mE\int^{t\wedge\tau}_0
\|\sigma^n_s(X_s)-\sigma_s(X_s)\|^2\dif s\\
&\leq\varliminf_{\eps\downarrow 0}
\frac{1}{\eps^2}\lim_{n\to\infty}\||\sigma^n-\sigma|^2\|_{\mL^{q/2}_{p/2}(t)}
=\varliminf_{\eps\downarrow 0}
\frac{1}{\eps^2}\lim_{n\to\infty}\|\sigma^n-\sigma\|^2_{\mL^q_p(t)}=0,
\end{align*}
and also,
$$
I_3(t)=0.
$$
For $I_1(t)$, we have
\begin{align*}
I_1(t)&\stackrel{(\ref{Es2})}{\leq} C\sup_{n\in\mN}\mE\int^{t\wedge\tau}_0
\Big[\cM|\nabla\sigma^n_s|(X_s)+\cM|\nabla\sigma^n_s|(Y_s)\Big]^2\dif s\\
&\leq C\sup_{n\in\mN}\|(\cM|\nabla\sigma^n_\cdot|)^2\|_{\mL^{q/2}_{p/2}(t)}
=C\sup_{n\in\mN}\|\cM|\nabla\sigma^n_\cdot|\|^2_{\mL^q_p(t)}\\
&\stackrel{(\ref{Es30})}{\leq} C\sup_{n\in\mN}\|\nabla\sigma^n_\cdot\|^2_{\mL^q_p(t)}\leq C\|\nabla\sigma_\cdot\|^2_{\mL^q_p(t)}.
\end{align*}
Combining the above calculations, we obtain that for all $t\geq 0$,
\begin{align}
\mE A_{t\wedge\tau}\leq C\|\nabla\sigma_\cdot\|^2_{\mL^q_p(t)}.\label{BP2}
\end{align}
Similarly, we can prove that
$$
\mE|M_{t\wedge\tau}|^2=4\mE\int^{t\wedge\tau}_0\frac{|[\sigma_s(X_s)-\sigma_s(Y_s)]^*Z_s|^2}{|Z_s|^4}\dif s\leq
C\|\nabla\sigma_\cdot\|^2_{\mL^q_p(t)},
$$
where the star denotes the transpose of a matrix.
The existence of a unique strong solution now follows from the classical Yamada-Watanabe theorem (cf. \cite{Ik-Wa}).
\end{proof}

Below, we prove better regularities of solutions with respect to the initial values.
\bl\label{Le5}
Under {\bf (H$^\sigma_1$)} and {\bf (H$^\sigma_2$)}, let $X_t(x)$ be the unique strong solution of SDE (\ref{DSDE}).
For any $T>0$, $\gamma\in\mR$ and all $x\not= y\in\mR^d$, we have
$$
\sup_{t\in[0,T]}\mE\left(|X_t(x)-X_t(y)|^{2\gamma}\right)\leq C|x-y|^{2\gamma},
$$
where $C=C(K,\delta,p,q,d,\gamma,T)$.
\el
\begin{proof}

For $x\not=y$ and $\eps\in(0,|x-y|)$, define
$$
\tau_\eps:=\inf\{t\geq 0: |X_t(x)-X_t(y)|\leq\eps\}.
$$
Set $Z^\eps_t:=X_{t\wedge\tau_\eps}(x)-X_{t\wedge\tau_\eps}(y)$.
For any $\gamma\in\mR$, by It\^o's formula, we have
\begin{align*}
|Z^\eps_t|^{2\gamma}&=|x-y|^{2\gamma}+2\gamma\int^{t\wedge\tau_\eps}_0|Z^\eps_s|^{2(\gamma-1)}\<Z^\eps_s,[\sigma_s(X_s(x))-\sigma_s(X_s(y))]\dif W_s\>\\
&\quad+2\gamma\int^{t\wedge\tau_\eps}_0|Z^\eps_s|^{2(\gamma-1)}\|\sigma_s(X_s(x))-\sigma_s(X_s(y))\|^2\dif s\\
&\quad+2\gamma(\gamma-1)\int^{t\wedge\tau_\eps}_0|Z^\eps_s|^{2(\gamma-2)}|[\sigma_s(X_s(x))-\sigma_s(X_s(y))]^*Z^\eps_s|^2\dif s\\
&=:|x-y|^{2\gamma}+\int^{t\wedge\tau_\eps}_0|Z^\eps_s|^{2\gamma}\Big(\alpha(s)\dif W_s+\beta(s)\dif s\Big),
\end{align*}
where
$$
\alpha(s):=\frac{2\gamma[\sigma_s(X_s(x))-\sigma_s(Y_s(y))]^*Z^\eps_s}{|Z^\eps_s|^2}
$$
and
$$
\beta(s):=\frac{2\gamma\|\sigma_s(X_s(x))-\sigma_s(Y_s(y))\|^2}{|Z^\eps_s|^2}
+\frac{2\gamma(\gamma-1)|[\sigma_s(X_s(x))-\sigma_s(Y_s(y))]^*Z^\eps_s|^2}{|Z^\eps_s|^4}.
$$
By the Dol\'eans-Dade's exponential (cf. \cite{Pr}), we have
$$
|Z^\eps_t|^{2\gamma}=|x-y|^{2\gamma}\exp\left\{\int^{t\wedge\tau_\eps}_0\alpha(s)\dif W_t
-\frac{1}{2}\int^{t\wedge\tau_\eps}_0|\alpha(s)|^2\dif s+\int^{t\wedge\tau_\eps}_0\beta(s)\dif s\right\}.
$$
Fix $T>0$ below.
Using (\ref{Es1}) and as in the proof of (\ref{BP2}), we have for any $0\leq s<t\leq T$,
$$
\mE\left(\int^t_s|\beta(r\wedge\tau_\eps)|\dif r\Bigg|_{\sF_s}\right)\leq  C\|\nabla\sigma\|^2_{\mL^{q}_{p}(s,t)},
$$
where $C=C(K,\delta,p,q,d,\gamma,T)$. Thus, by Lemma \ref{Le1}, we get for any $\lambda>0$,
$$
\mE\exp\left(\lambda\int^{T\wedge\tau_\eps}_0|\beta(s)|\dif s\right)\leq
\mE\exp\left(\lambda\int^T_0|\beta(s\wedge\tau_\eps)|\dif s\right)<+\infty.
$$
Similarly, we have
$$
\mE\exp\left(\lambda\int^{T\wedge\tau_\eps}_0|\alpha(s)|^2\dif s\right)<+\infty,\ \ \forall\lambda>0.
$$
In particular, by Novikov's criterion,
$$
t\mapsto\exp\left\{2\int^{t\wedge\tau_\eps}_0\alpha(s)\dif W_s
-2\int^{t\wedge\tau_\eps}_0|\alpha(s)|^2\dif s\right\}=:M^\eps_t
$$
is a continuous exponential martingale.
Hence, by H\"older's inequality, we have
$$
\mE|Z^\eps_t|^{2\gamma}\leq|x-y|^{2\gamma}(\mE M^\eps_t)^{\frac{1}{2}}
\left(\mE\exp\left\{\int^{t\wedge\tau_\eps}_0|\alpha(s)|^2\dif s
+2\int^{t\wedge\tau_\eps}_0\beta(s)\dif s\right\}\right)^{\frac{1}{2}}
\leq C|x-y|^{2\gamma},
$$
where $C$ is independent of $\eps$ and $x,y$.

Noting that
$$
\lim_{\eps\downarrow 0}\tau_\eps=\tau:=\inf\{t\geq 0: X_t(x)=X_t(y)\},
$$
by Fatou's lemma, we obtain
$$
\mE|X_{t\wedge\tau}(x)-X_{t\wedge\tau}(y)|^{2\gamma}=\varliminf_{\eps\to 0}\mE|Z^\eps_t|^{2\gamma}\leq C|x-y|^{2\gamma}.
$$
Letting $\gamma=-1$ yields that
$$
\tau\geq t,\ \ a.s.
$$
The proof is thus complete.
\end{proof}
Since $\sigma$ is bounded, the following lemma is standard, and we omit the details.
\bl\label{Le6}
Under {\bf (H$^\sigma_1$)}, let $X_t(x)$ solve SDE (\ref{DSDE}).
For any $T>0$, $\gamma\in\mR$ and all $x\in\mR^d$, we have
$$
\mE\left(\sup_{t\in[0,T]}(1+|X_t(x)|^2)^\gamma\right)\leq C_1(1+|x|^2)^\gamma,
$$
where $C_1=C_1(K,\gamma,T)$, and for any $\gamma\geq 1$ and $t,s\geq 0$,
$$
\sup_{x\in\mR^d}\mE|X_t(x)-X_s(x)|^{2\gamma}\leq C_2|t-s|^\gamma,
$$
where $C_2=C_2(K,\gamma)$.
\el

Basing on Lemmas \ref{Le5} and \ref{Le6}, it is  by now standard to prove the following theorem (cf. \cite[Theorem 4.5.1]{Ku}).
For the reader's convenience, we sketch the proof here.
\bt\label{Dif}
Under {\bf (H$^\sigma_1$)} and {\bf (H$^\sigma_2$)}, let $X_t(x)\in\sS^\infty_{0,\sigma}(x)$
be the unique strong solution of SDE (\ref{DSDE}),
then for almost all $\omega$ and all $t\in\mR_+$,
$x\mapsto X_t(\omega,x)$ is a homeomorphism on $\mR^d$.
\et
\begin{proof}
For $x\not=y\in\mR^d$, define
$$
\cR_t(x,y):=|X_{t}(x)-X_{t}(y)|^{-1}.
$$
For any  $x,y, x',y'\in\mR^d$ with $x\not= y$, $x'\not=y'$ and $s\not= t$,
it is easy to see that
$$
|\cR_t(x,y)-\cR_s(x',y')|\leq \cR_t(x,y)\cdot \cR_s(x',y')
\cdot[|X_t(x)-X_s(x')|+|X_t(y)-X_s(y')|].
$$
By Lemmas \ref{Le5} and \ref{Le6}, for any $\gamma\geq 1$ and $s,t\in[0,T]$, we have
$$
\mE|\cR_t(x,y)-\cR_s(x',y')|^\gamma\leq C|x-y|^{-\gamma}|x'-y'|^{-\gamma}(|t-s|^{\gamma/2}+|x-x'|^\gamma+|y-y'|^\gamma).
$$
Choosing $\gamma>4(d+1)$, by Kolmogorov's continuity criterion, there exists a continuous version to the mapping
$(t,x,y)\mapsto\cR_t(x,y)$ on $\{(t,x,y)\in\mR_+\times\mR^d\times\mR^d: x\not=y\}$. In particular, this proves that for almost
all $\omega$, the mapping $x\mapsto X_t(\omega,x)$ is one-to-one for all $t\geq 0$.

As for the onto property, let us define
$$
\cJ_t(x)=\left\{
\begin{array}{ll}
(1+|X_t(x|x|^{-2})|)^{-1},&\ \ x\not=0,\\
0,&\ \ x=0.
\end{array}
\right.
$$
As above, using Lemmas \ref{Le5} and \ref{Le6}, one can show that $(t,x)\mapsto\cJ_t(x)$ admits a continuous version.
Thus, $(t,x)\mapsto X_t(\omega,x)$ can be extended to a continuous map from $\mR_+\times\hat\mR^d$ to $\hat\mR^d$, where
$\hat\mR^d=\mR^d\cup\{\infty\}$ is the one-point compactification of $\mR^d$. Hence,
$X_t(\omega,\cdot):\hat\mR^d\to\hat\mR^d$ is homotopic to the identity mapping $X_0(\cdot)$ so that it is an onto map by
the well known fact in homotopic theory. In particular, for almost all $\omega$,
$x\mapsto X_t(\omega,x)$ is a homeomorphism on $\hat\mR^d$ for all $t\geq 0$.
Clearly, the restriction of $X_t(\omega,\cdot)$ to $\mR^d$ is still a homeomorphism
since $X_t(\omega,\infty)=\infty$.
\end{proof}

Now we turn to the proof of the strong Feller property.
\bt\label{Fel}
Under {\bf (H$^\sigma_1$)} and {\bf (H$^\sigma_2$)}, let $X_t(x)\in\sS^\infty_{0,\sigma}(x)$
be the unique strong solution of SDE (\ref{DSDE}), then for any bounded measurable function $\phi$, $T>0$ and $x,y\in\mR^d$,
\begin{align}
|\mE(\phi(X_t(x)))-\mE(\phi(X_t(y)))|\leq \frac{C_T}{\sqrt{t}}\|\phi\|_\infty|x-y|,\ \ \forall t\in(0,T].\label{BB10}
\end{align}
\et
\begin{proof}
Define $\sigma^n_t(x):=\sigma_t*\rho_n(x)$, where $\rho_n$ is a mollifier in $\mR^d$. By {\bf (H$^\sigma_1$)},
it is easy to see that for all $(t,x)\in\mR_+\times\mR^d$,
\begin{align}
\delta|\lambda|^2\leq \sum_{ik}|[\sigma^n_t(x)]^{ik}\lambda^i|^2\leq K|\lambda|^2,\ \ \forall\lambda\in\mR^d.\label{MM3}
\end{align}
Let $X^n_t(x)\in\sS^\infty_{0,\sigma^n}(x)$ be the unique strong solution of SDE (\ref{DSDE})
corresponding to $\sigma^n$. By the monotone class theorem, it suffices to prove (\ref{BB10}) for any bounded Lipschitz
continuous function $\phi$.
First of all, by Bismut-Elworthy-Li's formula (cf. \cite{El-Li}), for any $h\in\mR^d$, we have
\begin{align}
\nabla_h\mE\phi(X^n_t(x))=\frac{1}{t}\mE\left[\phi(X^n_t(x))\int^t_0[\sigma^n_s(X^n_s(x))]^{-1}\nabla_h X^n_s(x)\dif W_s\right],\label{MM1}
\end{align}
where for a smooth function $f$, we denote $\nabla_hf:=\<\nabla f, h\>$.
Noting that
$$
\nabla_h X^n_t(x)=h+\int^t_0\nabla\sigma^n_s(X^n_s(x))\cdot\nabla_h X^n_s(x)\dif W_s,
$$
by It\^o's formula, we have
\begin{align*}
|\nabla_h X^n_t(x)|^2&=|h|^2+2\int^t_0\<\nabla_h X^n_s(x),\nabla\sigma^n_s(X^n_s(x))\cdot\nabla_h X^n_s(x)\dif W_s\>\\
&\quad+\int^t_0\|\nabla\sigma^n_s(X^n_s(x))\cdot\nabla_h X^n_s(x)\|^2\dif s\\
&=:|h|^2+\int^t_0|\nabla_h X^n_s(x)|^2\Big(\a^n_h(s)\dif W_s+\beta^n_h(s)\dif s\Big),
\end{align*}
where
$$
\alpha^n_h(s):=\frac{(\nabla_h X^n_s(x))^*\cdot\nabla\sigma^n_s(X^n_s(x))\cdot\nabla_h X^n_s(x)}{|\nabla_h X^n_s(x)|^2}
$$
and
$$
\beta^n_h(s):=\frac{\|\nabla\sigma^n_s(X^n_s(x))\cdot\nabla_h X^n_s(x)\|^2}{|\nabla_h X^n_s(x)|^2}.
$$
By the Dol\'eans-Dade's exponential again, we have
$$
|\nabla_h X^n_t(x)|^2=|h|^2\exp\left\{\int^t_0\alpha^n_h(s)\dif W_s
-\frac{1}{2}\int^t_0|\alpha^n_h(s)|^2\dif s+\int^t_0\beta^n_h(s)\dif s\right\}.
$$
Fix $T>0$. By (\ref{Es1}), we have for any $0\leq s<t\leq T$,
$$
\mE\left(\int^t_s|\beta^n_h(r)|\dif r\Bigg|_{\sF_s}\right)\leq  C\|\nabla\sigma^n\|^2_{\mL^{q}_{p}(s,t)}\leq
C\|\nabla\sigma\|^2_{\mL^{q}_{p}(s,t)},
$$
where $C=C(K,\delta,p,q,d,T)$ is independent of $n,x$ and $h$. Thus, by Lemma \ref{Le1}, we get for any $\lambda>0$,
$$
\sup_n\sup_{h\in\mR^d}\mE\exp\left(\lambda\int^T_0|\beta^n_h(s)|\dif s\right)<+\infty.
$$
Similarly,
$$
\sup_n\sup_{h\in\mR^d}\mE\exp\left(\lambda\int^T_0|\alpha^n_h(s)|^2\dif s\right)<+\infty.
$$
Hence,
$$
\sup_n\sup_{t\in[0,T]}\sup_{x\in\mR^d}\mE|\nabla_h X^n_t(x)|^2\leq C|h|^2,\ \ \forall h\in\mR^d,
$$
and by (\ref{MM3}) and (\ref{MM1}),
\begin{align*}
|\nabla_h\mE\phi(X^n_t(x))|&\leq\frac{\|\phi\|_\infty}{t}\left(\mE\int^t_0|[\sigma^n_s(X^n_s(x))]^{-1}\nabla_h X^n_s(x)|^2\dif s\right)^{\frac{1}{2}}\\
&\leq\frac{C_T\|\phi\|_\infty}{t}\left(\mE\int^t_0|\nabla_h X^n_s(x)|^2\dif s\right)^{\frac{1}{2}}
\leq\frac{C_T\|\phi\|_\infty|h|}{\sqrt{t}},
\end{align*}
which implies that for all $t\in(0,T]$ and $x,y\in\mR^d$,
\begin{align}
|\mE(\phi(X^n_t(x)))-\mE(\phi(X^n_t(y)))|\leq \frac{C_T\|\phi\|_\infty}{\sqrt{t}}|x-y|,\label{BB101}
\end{align}
where $C_T$ is independent of $n$.

Now for completing the proof, it only needs to take limits for (\ref{BB101}) by proving that for any $x\in\mR^d$,
\begin{align}
\lim_{n\to\infty}\mE|X^n_t(x)-X_t(x)|=0.\label{MM2}
\end{align}
Set
$$
Z^n_t(x):=X^n_t(x)-X_t(x)
$$
and
$$
\eta^n(s):=\Big(\cM|\nabla\sigma^n_s|(X^n_s(x))+\cM|\nabla\sigma^n_s|(X_s(x))\Big)^2.
$$
For any $\lambda>0$, by It\^o's formula, we have
\begin{align*}
\mE|Z^n_t(x)|^2 e^{-\lambda\int^t_0\eta^n(s)\dif s}&=
\mE\int^t_0\|\sigma^n_s(X^n_s(x))-\sigma_s(X_s(x))\|^2e^{-\lambda\int^s_0\eta^n(r)\dif r}\dif s\\
&\quad-\lambda\mE\int^t_0\eta^n(s)|Z^n_s(x)|^2 e^{-\lambda\int^s_0\eta^n(r)\dif r}\dif s\\
&\leq\mE\int^t_0\|\sigma^n_s(X^n_s(x))-\sigma^n_s(X_s(x))\|^2e^{-\lambda\int^s_0\eta^n(r)\dif r}\dif s\\
&\quad+\mE\int^t_0\|\sigma^n_s(X_s(x))-\sigma_s(X_s(x))\|^2e^{-\lambda\int^s_0\eta^n(r)\dif r}\dif s\\
&\quad-\lambda\mE\int^t_0\eta^n(s)|Z^n_s(x)|^2 e^{-\lambda\int^s_0\eta^n(r)\dif r}\dif s\\
&\stackrel{(\ref{Es2})}{\leq} (C_d-\lambda)\mE\int^t_0\eta^n(s)|Z^n_s(x)|^2 e^{-\lambda\int^s_0\eta^n(r)\dif r}\dif s\\
&\quad+\mE\int^t_0\|\sigma^n_s(X_s(x))-\sigma_s(X_s(x))\|^2\dif s.
\end{align*}
Thus, by (\ref{Es1}), we obtain that for any $\lambda\geq C_d$,
$$
\lim_{n\to\infty}\mE|Z^n_t(x)|^2 e^{-\lambda\int^t_0\eta^n(s)\dif s}\leq\lim_{n\to\infty}
\|\sigma^n-\sigma\|^2_{\mL^q_p(T)}=0.
$$
Moreover, as above, by (\ref{Es1}), (\ref{Es30}) and Lemma \ref{Le1}, we also have
$$
\sup_n\mE\exp\left(\lambda\int^T_0|\eta^n(s)|\dif s\right)<+\infty,\ \ \forall\lambda,T>0.
$$
Hence, by H\"older's inequality,
$$
\lim_{n\to\infty}\mE|Z^n_t(x)|\leq \lim_{n\to\infty}\left[\left(\mE e^{\lambda\int^t_0\eta^n(s)\dif s}\right)^{\frac{1}{2}}
\left(\mE|Z^n_t(x)|^2 e^{-\lambda\int^t_0\eta^n(s)\dif s}\right)^{\frac{1}{2}}\right]=0,
$$
which then gives (\ref{MM2}). The proof is complete.
\end{proof}

\section{Zvonkin's transformation and Proof of Theorem \ref{Main}}

In this section we prove Theorem \ref{Main} by using Zvonkin's transformation to kill the drift (cf. \cite{Zv}).
Below, we assume that $\sigma$ satisfies {\bf (H$^\sigma_1$)} and $b\in L^q(\mR_+,L^p(\mR^d))$ provided with
\begin{align}
\frac{d}{p}+\frac{2}{q}<1.\label{BB88}
\end{align}

Fix $T_0>0$. For any $T\in[0,T_0]$ and $\ell=1,\cdots,d$, let $u^\ell(t,x)$ solve the following PDE:
$$
\p_tu^\ell(t,x)+L_tu^\ell(t,x)+b^\ell(t,x)=0,\ \ u^\ell(T,x)=0,
$$
where $L_t$ is given by (\ref{LL2}).
Set
$$
\u(t,x):=(u^1(t,x),\cdots,u^d(t,x))\in\mR^d.
$$
By Theorem \ref{Th11}, we have
\begin{align}
C_0:=\sup_{T\in[0,T_0]}\Big(\|\p_t\u\|_{\mL^q_p(T)}+\|\u\|_{\mH^{2,q}_p(T)}\Big)<+\infty.\label{BP8}
\end{align}
Thanks to (\ref{BB88}) and (\ref{BP8}),
by \cite[Lemma 10.2]{Kr-Ro},
$$
(t,x)\mapsto\nabla\u(t,x)\mbox{ is H\"older continuous},
$$
and for fixed $\delta\in(0,\frac{1}{2}-\frac{d}{2p}-\frac{1}{q})$,
there exists constant $C_1>0$ depending only on $p,q,\delta$ such that for any $S\in[0,T]$,
\begin{align}
\sup_{(t,x)\in[S,T]\times\mR^d}|\nabla\u(t,x)|\leq C_1(T-S)^\delta
\Big(\|\p_t\u\|_{\mL^q_p(S,T)}+\|\u\|_{\mH^{2,q}_p(S,T)}\Big)\leq C_0C_1 (T-S)^\delta,\label{LL3}
\end{align}
where $C_0$ is defined by (\ref{BP8}).

Let $\u_n$ be the mollifying approximation of $\u$ defined as in (\ref{BP6}). Define
$$
\Phi_t(x):=x+\u(t,x),\ \ \ \Phi^n_t(x):=x+\u_n(t,x).
$$
It is easy to see that $\Phi$ solves the following PDE:
\begin{align}
\p_t\Phi_t(x)+L_t\Phi_t(x)=0,\ \ \Phi_T(x)=x.\label{BB1}
\end{align}
Moreover, letting $T,S\in[0, T_0]$ satisfy that
\begin{align}
0\leq T-S\leq \frac{1}{2(C_0C_1)^{1/\delta}},\label{LL4}
\end{align}
then by (\ref{LL3}), we have for all $t\in[S,T]$,
$$
\tfrac{1}{2}|x-y|\leq|\Phi^n_t(x)-\Phi^n_t(y)|\leq \tfrac{3}{2}|x-y|
$$
and
$$
\tfrac{1}{2}|x-y|\leq|\Phi_t(x)-\Phi_t(y)|\leq \tfrac{3}{2}|x-y|,
$$
which implies that $\Phi_t$ and $\Phi^n_t$ are diffeomorphisms on $\mR^d$. So, if we set
$$
\Psi_t(x):=\Phi^{-1}_t(x),\ \ \Psi^n_t(x):=\Phi^{n,-1}_t(x),
$$
then
\begin{align}
|\nabla\Phi_t(x)|\vee|\nabla\Phi^n_t(x)|\leq\frac{3}{2},\ \ |\nabla\Psi_t(x)|\vee|\nabla\Psi^n_t(x)|\leq 2.\label{BP9}
\end{align}

We first prove two lemmas.
\bl\label{Le3}
For each $(t,x)\in[S,T]\times\mR^d$, we have
\begin{align}
\lim_{n\to\infty}\Phi^n_t(x)=\Phi_t(x),\ \ \lim_{n\to\infty}\Psi^n_t(x)=\Psi_t(x)\label{BP99}
\end{align}
and
\begin{align}
\lim_{n\to\infty}|\nabla\Psi^n_t(y)-\nabla\Psi_t(y)|=0.\label{BP10}
\end{align}
\el
\begin{proof}
The first limit is immediate from the property of convolution. The second limit follows from
$$
|\Psi^n_t(x)-\Psi_t(x)|\leq 2|x-\Phi^n_t(\Psi_t(x))|=2|\Phi_t(\Psi_t(x))-\Phi^n_t(\Psi_t(x))|,
$$
and the first limit. As for the third limit, noting that
$$
[\nabla\Psi^n_t(y)]^{-1}=\nabla\Phi^n_t\circ\Psi^n_t(y),
$$
by (\ref{BP9}), we have
\begin{align*}
|\nabla\Psi^n_t(y)-\nabla\Psi_t(y)|&=|\nabla\Psi^n_t(y)|\cdot|\nabla\Phi^n_t\circ\Psi^n_t(y)-
\nabla\Phi_t\circ\Psi_t(y)|\cdot|\nabla\Psi_t(y)|\\
&\leq4|\nabla\Phi^n_t\circ\Psi^n_t(y)-\nabla\Phi_t\circ\Psi_t(y)|.
\end{align*}
The third limit follows from the continuity of $x\mapsto\nabla\Phi_t(x)$.
\end{proof}
\bl\label{Le4}
We have
$$
\lim_{n\to\infty}\|\p_i\Psi^{n,i'}_s\cdot\p_j\Psi^{n,j'}_s\cdot(\p_{i'}\p_{j'}\Phi^{n,l}_s\circ\Psi^n_s)
\cdot\p_l\Psi^{n,k}_s-\p_i\Psi^{i'}_s\cdot\p_j\Psi^{j'}_s\cdot(\p_{i'}\p_{j'}\Phi^{l}_s\circ\Psi_s)
\cdot\p_l\Psi^{k}_s\|_{\mL^q_p(S,T)}=0
$$
and
$$
\lim_{n\to\infty}\|(\p_t\Phi^n\circ\Psi^n)\cdot\nabla\Psi^n-(\p_t\Phi\circ\Psi)\cdot\nabla\Psi\|_{\mL^q_p(S,T)}=0.
$$
\el
\begin{proof}
We only prove the first limit, the second limit can be proved similarly. For proving the first limit, it suffices to prove the following two limits:
$$
\lim_{n\to\infty}\|\p_i\Psi^{n,i'}_s\cdot\p_j\Psi^{n,j'}_s\cdot\p_{i'}\p_{j'}\Phi^{l}_s\circ\Psi_s
\cdot\p_l\Psi^{n,k}_s-\p_i\Psi^{i'}_s\cdot\p_j\Psi^{j'}_s\cdot\p_{i'}\p_{j'}\Phi^{l}_s\circ\Psi_s
\cdot\p_l\Psi^{k}_s\|_{\mL^q_p(S,T)}=0,
$$
$$
\lim_{n\to\infty}\|\p_i\Psi^{n,i'}_s\cdot\p_j\Psi^{n,j'}_s\cdot\p_{i'}\p_{j'}\Phi^{n,l}_s\circ\Psi^n_s
\cdot\p_l\Psi^{n,k}_s-\p_i\Psi^{n,i'}_s\cdot\p_j\Psi^{n,j'}_s\cdot\p_{i'}\p_{j'}\Phi^{l}_s\circ\Psi_s
\cdot\p_l\Psi^{n,k}_s\|_{\mL^q_p(S,T)}=0.
$$
The first limit follows by (\ref{BP8}), (\ref{BP9}), (\ref{BP10}) and the dominated convergence theorem.
For the second limit, by (\ref{BP9}), we have
\begin{align*}
&\|\p_i\Psi^{n,i'}_s\cdot\p_j\Psi^{n,j'}_s\cdot\p_{i'}\p_{j'}\Phi^{n,l}_s\circ\Psi^n_s
\cdot\p_l\Psi^{n,k}_s-\p_i\Psi^{n,i'}_s\cdot\p_j\Psi^{n,j'}_s\cdot\p_{i'}\p_{j'}\Phi^{l}_s\circ\Psi_s
\cdot\p_l\Psi^{n,k}_s\|_{\mL^q_p(S,T)}\\
&\qquad\leq
8\|\nabla^2\Phi^n\circ\Psi^n-\nabla^2\Phi\circ\Psi\|_{\mL^q_p(S,T)}=
8\|\nabla^2\u_n\circ\Psi^n-\nabla^2\u\circ\Psi\|_{\mL^q_p(S,T)}\leq\\
&\qquad\leq 8\|\nabla^2\u_n\circ\Psi^n-\nabla^2\u\circ\Psi^n\|_{\mL^q_p(S,T)}
+8\|\nabla^2\u\circ\Psi^n-\nabla^2\u\circ\Psi\|_{\mL^q_p(S,T)}\\
&\qquad\leq C\|\nabla^2\u_n-\nabla^2\u\|_{\mL^q_p(S,T)}
+8\|\nabla^2\u\circ\Psi^n-\nabla^2\u\circ\Psi\|_{\mL^q_p(S,T)}.
\end{align*}
where in the last step, we have used the change of variables and
$$
\sup_n\sup_{(t,x)\in[S,T]\times\mR^d}\det(\nabla\Phi^n_t(x))\leq C.
$$
It is clear that by (\ref{BP8}),
$$
\lim_{n\to\infty}\|\nabla^2\u_n-\nabla^2\u\|_{\mL^q_p(S,T)}=0.
$$
On the other hand, let $\u_\eps$ be a family of smooth functions on $[0,T]\times\mR^d$ with compact supports such that
$$
\lim_{\eps\to 0}\|\nabla^2\u_\eps-\nabla^2\u\|_{\mL^q_p(S,T)}=0.
$$
Then as above, we have
$$
\lim_{\eps\to 0}\sup_n\|\nabla^2\u_\eps\circ\Psi^n-\nabla^2\u\circ\Psi^n\|_{\mL^q_p(S,T)}=0,
$$
and for fixed $\eps$, by (\ref{BP99}) and the dominated convergence theorem,
$$
\lim_{n\to\infty}\|\nabla^2\u_\eps\circ\Psi^n-\nabla^2\u_\eps\circ\Psi\|_{\mL^q_p(S,T)}=0.
$$
Hence,
$$
\lim_{n\to\infty}\|\nabla^2\u\circ\Psi^n-\nabla^2\u\circ\Psi\|_{\mL^q_p(S,T)}=0.
$$
The proof is thus complete.
\end{proof}
Now we are in a position to prove the following Zvonkin's transformation to kill the drift.
\bl\label{Le7}
Let $\tau$ be any ($\sF_t$)-stopping time.
Let $X_t$ be a $\mR^d$-valued ($\sF_t$)-adapted and continuous stochastic process satisfying
$$
P\left\{\omega: \int^t_0\Big(|b_s(X_s(\omega))|+|\sigma_s(X_s(\omega))|^2\Big)\dif s<+\infty,
\forall t\in[0,\tau(\omega))\right\}=1.
$$
Then $X_t$ solves the following SDE on $[S\wedge\tau,T\wedge\tau)$,
$$
\dif X_t=b_t(X_t)\dif t+\sigma_t(X_t)\dif W_t,
$$
if and only if $Y_t:=\Phi_t(X_t)$ solves the following SDE on $[S\wedge\tau,T\wedge\tau)$
$$
\dif Y_t=\Sigma_t(Y_t)\dif W_t,
$$
where $\Sigma^{ik}_t(y):=(\p_l\Phi^i_t\cdot\sigma^{lk}_t)\circ\Psi_t(y)$.
\el
\begin{proof}
We first prove the ``only if'' part. Let $X^n_t:=\Psi^n_t(Y_t).$
By It\^o's formula, we have for all $t\in[S\wedge\tau,T\wedge\tau)$,
\begin{align}
X^n_t=\Psi^n_{S\wedge\tau}(Y_{S\wedge\tau})+\int^t_{S\wedge\tau}\Big[\p_s\Psi^n_s
+\tfrac{1}{2}(\Sigma_s\Sigma^*_s)^{ij}\p_i\p_j\Psi^n_s\Big](Y_s)\dif s
+\int^t_{S\wedge\tau}[\nabla\Psi^n_s\cdot\Sigma_s](Y_s)\dif W_s.\label{St}
\end{align}
Noticing that
$$
\p_s\Psi^{n}_s\cdot(\nabla\Phi^n_s\circ\Psi^n_s)+\p_s\Phi^n_s\circ\Psi^n_s=0
$$
and
$$
\p_i\Psi^{n,i'}_s\cdot\p_j\Psi^{n,j'}_s\cdot(\p_{i'}\p_{j'}\Phi^{n,l}_s\circ\Psi^n_s)
+\p_i\p_j\Psi^{n,k}_s\cdot(\p_k\Phi^{n,l}_s\circ\Psi^n_s)=0,
$$
we have
$$
\p_s\Psi^{n}_s=-(\p_s\Phi^n_s\circ\Psi^n_s)\cdot \nabla\Psi^n_s
$$
and
$$
\p_i\p_j\Psi^{n,k}_s=-\p_i\Psi^{n,i'}_s\cdot\p_j\Psi^{n,j'}_s\cdot(\p_{i'}\p_{j'}\Phi^{n,l}_s\circ\Psi^n_s)
\cdot\p_l\Psi^{n,k}_s.
$$
Let $X_t=\Psi_t(Y_t)$.
Taking limits for both sides of (\ref{St}),
and by Lemmas \ref{Le3}, \ref{Le4} and (\ref{Es1}), (\ref{BB1}), one finds that  for all $t\in[S\wedge\tau,T\wedge\tau)$,
$$
X_t=\Psi_S(Y_{S\wedge\tau})+\int^t_{S\wedge\tau}b(X_s)\dif s
+\int^t_{S\wedge\tau}\sigma_s(X_s)\dif W_s.
$$
The ``if'' part is similar by (\ref{BB9}) and in fact easier. We omit the details.
\end{proof}

Basing on the above Zvonkin's transformation, we can give

\begin{proof}[Proof of Theorem \ref{Main}]
Using the standard time shift technique (cf. \cite{Zh2}), by Lemma \ref{Le7} and Theorems \ref{Dif}, \ref{Fel},
it only needs to check that $\Sigma^{ik}_t(y):=(\p_l\Phi^i_t\cdot\sigma^{lk}_t)\circ\Psi_t(y)$
 satisfies {\bf (H$^\Sigma_1$)} and {\bf (H$^\Sigma_2$)}. First of all,
{\bf (H$^\Sigma_1$)} is clear. For {\bf (H$^\Sigma_2$)}, we have
$$
\p_l\Sigma^{ik}_t(y)=[(\p_{l'}\p_j\Phi^i_t\cdot\sigma^{jk}_t+\p_j\Phi^i_t\cdot\p_{l'}\sigma^{jk}_t)\circ\Psi_t(y)]\cdot\p_l\Psi^{l'}(y).
$$
By (\ref{BP8}), (\ref{BP9}) and  {\bf (H$^\sigma_2$)}, it is easy to see that
$$
\|\p_l\Sigma^{ik}\|_{\mL^q_p(T_0)}<+\infty.
$$
\end{proof}

\section{Appendix}

The following result is a combination of \cite[Theorem 10.3 and Remark 10.4]{Kr-Ro}.
\bt\label{Th11}
Let $p,q\in(1,\infty)$ satisfy (\ref{BP4}). Assume {\bf (H$^\sigma_1$)} and $b\in L^q(\mR_+,L^p(\mR^d))$.
For any $T>0$ and $f\in \mL^q_p(T)$, there exists a unique solution
$u\in \cH^{2,q}_p(T)$ for the following PDE:
\begin{align}
\p_tu(t,x)+L_tu(t,x)+f(t,x)=0,\ \ u(T,x)=0.\label{PDE}
\end{align}
Moreover, this solution satisfies that for any $S\in[0,T]$,
$$
\|\p_tu\|_{\mL^q_p(S,T)}+\|u\|_{\mH^{2,q}_p(S,T)}\leq C\|f\|_{\mL^q_p(S,T)},
$$
where $C=C(T,K,\delta,p,q, \|b\|_{\mL^q_p(T)})$.
\et
The following result can be proved along the same lines as in \cite[Theorem 10.3, Remark 10.4]{Kr-Ro}. We omit the details.
\bt\label{Th1}
Assume {\bf (H$^\sigma_1$)} and we consider the following two cases about $b$:
\begin{enumerate}[{\bf (1$^o$)}]
\item Let $p,q\in(1,\infty)$ be fixed and let
$b$ be a bounded measurable function.
\item Let $p,q\in(1,\infty)$ satisfy (\ref{BP4}) and  let $b\in L^q(\mR_+,L^p(\mR^d))\cap L^\infty(\mR_+\times\mR^d)$.
\end{enumerate}
For any $T>0$, $r\in(1,\infty)$ and $f\in \mL^r_r(T)\cap \mL^q_p(T)$, there exists a unique solution
$u\in \cH^{2,r}_r(T)\cap\cH^{2,q}_p(T)$ for PDE (\ref{PDE}).
Moreover, this solution satisfies that for any $S\in[0,T]$,
$$
\|\p_tu\|_{\mL^r_r(S,T)}+\|u\|_{\mH^{2,r}_r(S,T)}\leq C_1\|f\|_{\mL^r_r(S,T)}
$$
and
$$
\|\p_tu\|_{\mL^q_p(S,T)}+\|u\|_{\mH^{2,q}_p(S,T)}\leq C_2\|f\|_{\mL^q_p(S,T)},
$$
where $C_1=C_1(T,K,\delta,p,q, \|b\|_\infty)$ and,
in case {\bf (1$^o$)}, $C_2=C_2(T,K,\delta,p,q, \|b\|_\infty)$, and in case {\bf (2$^o$)}, $C_2=C_2(T,K,\delta,p,q, \|b\|_{\mL^q_p(T)})$.
\et
The following lemma is taken from \cite[p. 1, Lemma 1.1]{Po}.
\bl\label{Le1}
Let $\{\beta(t)\}_{t\in[0,T]}$ be a nonnegative measurable ($\sF_t$)-adapted process. Assume that for all $0\leq s\leq t\leq T$,
$$
\mE\left(\int^t_s\beta(r)\dif r\Bigg|_{\sF_s}\right)\leq\rho(s,t),
$$
where $\rho(s,t)$ is a nonrandom interval function satisfying the following conditions:

(i) $\rho(t_1,t_2)\leq\rho(t_3,t_4)$ if $(t_1,t_2)\subset(t_3,t_4)$;

(ii) $\lim_{h\downarrow 0}\sup_{0\leq s<t\leq T, |t-s|\leq h}\rho(s,t)=\kappa,\ \ \kappa\geq0$.\\
Then for any arbitrary real $\lambda<\kappa^{-1}$ (if $\kappa=0$, then $\kappa^{-1}=+\infty$),
$$
\mE\exp\left\{\lambda\int^T_0\beta(r)\dif r\right\}\leq C=C(\lambda,\rho, T)<+\infty.
$$
\el
Let $\varphi$ be a locally integrable function
on $\mR^d$. The Hardy-Littlewood maximal function is defined by
$$
\cM \varphi(x):=\sup_{0<r<\infty}\frac{1}{|B_r|}\int_{B_r}\varphi(x+y)\dif y,
$$
where $B_r:=\{x\in\mR^d: |x|<r\}$.
The following result can be found in  \cite[Appendix A]{Cr-De-Le}.
\bl\label{Le2}
(i) There exists a constant $C_d>0$ such that for all $\varphi\in C^\infty(\mR^d)$ and $x,y\in \mR^d$,
\begin{align}
|\varphi(x)-\varphi(y)|\leq C_d\cdot |x-y|\cdot(\cM|\nabla\varphi|(x)+\cM|\nabla \varphi|(y)).\label{Es2}
\end{align}
(ii) For any $p>1$, there exists a constant $C_{d,p}$ such that for all $\varphi\in L^p(\mR^d)$,
\begin{align}
\left(\int_{\mR^d}(\cM\varphi(x))^p\dif x\right)^{1/p}
\leq C_{d,p}\left(\int_{\mR^d}|\varphi(x)|^p\dif x\right)^{1/p}.\label{Es30}
\end{align}
\el

\vspace{5mm}

{\bf Acknowledgements:}
The author is very grateful to Professor Michael R\"ockner for his valuable suggestions.
This work is supported by NSFs of China (No. 10971076; 10871215) and
Program for New Century Excellent Talents in University.

\end{document}